\documentclass[dvips,german,french,english,twoside]{article}
\usepackage[a5paper,layoutwidth=148mm,layoutheight=7.77817in,
outer=10mm,inner=22bp,tmargin=10mm,bmargin=10mm,footskip=15bp,
]{geometry}

\usepackage[utf8]{inputenc}
\usepackage[T1]{fontenc}
\usepackage{lmodern}
\usepackage{fix-cm}
\usepackage{textcomp}
\usepackage{textgreek}
\newcommand\textb[1]{{\fontseries{b}\selectfont #1}}
\newcommand\sub[2]{\textsl{#2}}
\newcommand\add[1]{\textsl{#1}}
\usepackage{natbib}
\bibpunct (),.{},

\newcommand{\mockalph}[1]{}
\makeatletter
\DeclareFontEncoding{LS1}{}{}
\DeclareFontEncoding{LS2}{}{\noaccents@}
\DeclareSymbolFont{stix-letters}       {LS1}{stix}     {m}{it}
\DeclareSymbolFont{stix-arrows1}       {LS1}{stixsf}   {m} {n}
\DeclareSymbolFont{stix-operators}     {LS1}{stix}     {m} {n}
\DeclareSymbolFont{stix-largesymbols}  {LS2}{stixex}   {m} {n}
\DeclareSymbolFont{stix-bold-operators}{LS1}{stix}     {b} {n}
\DeclareFontSubstitution{LS1}{stix}{m}{n}
\DeclareFontSubstitution{LS2}{stix}{m}{n}
\SetSymbolFont{stix-letters}     {bold}{LS1}{stix}     {b}{it}
\SetSymbolFont{stix-arrows1}     {bold}{LS1}{stixsf}   {b} {n}
\SetSymbolFont{stix-operators}   {bold}{LS1}{stix}     {b} {n}
\SetSymbolFont{stix-largesymbols}{bold}{LS2}{stixex}   {b} {n}

\def\stix@undefine#1{%
    \if\relax\noexpand#1\let#1=\@undefined\fi}
\def\stix@MathSymbol#1#2#3#4{%
    \stix@undefine#1%
    \DeclareMathSymbol{#1}{#2}{#3}{#4}}
\stix@MathSymbol{\stixwedge}{\mathbin}  {stix-operators}{"E1} \let\land=\wedge
\stix@MathSymbol{\stixvee}  {\mathbin}  {stix-operators}{"E2} \let\lor=\vee
\def\bigwedgegras{\DOTSI\bigwedgeop\slimits@}
\def\bigveegras{\DOTSI\bigveeop\slimits@}
\stix@MathSymbol{\stixbigwedgeop}{\mathop}{stix-largesymbols}{"B4}
\stix@MathSymbol{\stixbigveeop}  {\mathop}{stix-largesymbols}{"B5}
\stix@MathSymbol{\stixrightarrow}               {\mathrel}{stix-arrows1}{"99}
\stix@MathSymbol{\stixrightleftarrows}          {\mathrel}{stix-arrows1}{"CB}
\makeatother

\renewcommand\land{\mathbin{\boldsymbol{\stixwedge}}}
\renewcommand\lor{\mathbin{\boldsymbol{\stixvee}}}
\newcommand\bigland{\mathop{\textstyle
{\stixbigwedgeop}}\nolimits}
\newcommand\biglor{\mathop{\textstyle
{\stixbigveeop}}\nolimits}

\newcommand\Srel{\mathrel{S}}

\usepackage{babel} 
\newcommand\german[1]{\foreignlanguage{german}{#1}}
\usepackage[french]{varioref}
\usepackage{enumerate}
\usepackage{graphicx}
\usepackage{pdfpages}

\usepackage{xurl}
\usepackage[normalem]{ulem}

\usepackage{refcount}
\usepackage{ebproof}
\usepackage{amsmath}
\usepackage{amsfonts}
\usepackage{amssymb}

\usepackage{mathtools} \mathtoolsset{mathic} 
\newcommand*\dng[1]{%
  \setbox0\hbox{$\mathaccent"0362{#1}^H$}%
  \setbox2\hbox{$\mathaccent"0362{\kern0pt#1}^H$}%
  \ifdim\ht0=\ht2 \overline{\overline{#1}}\else \bar{\bar#1}\fi
  }

\usepackage{indentfirst}
\usepackage{framed}

\makeatletter
\newcommand{\leqnomode}{\tagsleft@true\let\veqno\@@leqno}
\newcommand{\reqnomode}{\tagsleft@false\let\veqno\@@eqno}
\makeatother

\newcommand\barre[2][]{#1}
\newcommand\bs{} 

\newcommand\Bs{-} 
\renewcommand\le\leqslant

\PassOptionsToPackage{urlcolor=magenta,citecolor=blue,anchorcolor=blue,linkcolor=blue,menucolor=cyan}{hyperref}%
\usepackage[pdftitle=,pdfauthor=Stefan\ Neuwirth,pdfdisplaydoctitle=true,pdflang=fr-FR,bookmarksopen=false,breaklinks=true,unicode=true,plainpages=false,hyperindex=true,pdfstartview=FitH,colorlinks=true,pdfpagelabels=true,bookmarksnumbered=true]{hyperref}
\ifx\texorpdfstring\undefined\newcommand\texorpdfstring[2]{#1}\fi
\ifx\pdfbookmark[3]{}\fi

\usepackage{doi}
\begin{document}
\title{In 1955, Paul Lorenzen clears the sky in foundations of mathematics for Hermann Weyl}
\author{Thierry Coquand \and Henri Lombardi \and Stefan Neuwirth}
\date{}
\maketitle


In 1955, Paul Lorenzen is a mathematician who devotes all his research to foundations of mathematics, on a par with Hans Hermes, but his academic background is algebra in the tradition of Helmut Hasse and Wolfgang Krull. This shift from algebra to logic goes along with his discovery that his ``algebraic works [\dots]\ have been concerned with a problem that has \emph{formally} the same structure as the problem of consistency of the classical calculus of logic'' (letter to \german{Carl Friedrich Gethmann} dated 14 January 1988, \citealp[see][page~76]{gethmann91}).

After having provided a proof of consistency for arithmetic in 1944 \citep[see][]{lorenzen20,coquandneuwirth19} and published it in 1951 \citep[see][]{lorenzen51wt,coquandneuwirth23}, Lorenzen inquires still further into the foundations of mathematics and arrives at the conviction that analysis can also be given a predicative foundation.

Wilhelm Ackermann as well as Paul Bernays have pointed out to him in 1947 that his views are very close to those proposed by Hermann Weyl in \emph{Das Kontinuum} (1918): sets are not postulated to exist beforehand; they are being generated in an ongoing process of comprehension.

This seems to be the reason for Lorenzen to get into contact with Weyl, who develops a genuine interest into Lorenzen's operative mathematics and welcomes with great enthusiasm his \emph{Einführung in die operative Logik und Mathematik} (1955), which he studies line by line. This book's aim is to grasp the objects of analysis by means of inductive definitions; the most famous achievement of this enterprise is a generalised inductive formulation of the Cantor-Bendixson theorem that makes it constructive.

This mathematical kinship is brutally interrupted by Weyl's death in 1955; a planned visit by Lorenzen at the Institute for Advanced Study in Princeton takes place only in 1957--1958.

As told by Kuno Lorenz, Lorenzen's first Ph.D. student, a discussion with Alfred Tarski during this visit provokes a turmoil in Lorenzen's operative research program that leads to his abandonment of language levels and to a great simplification of his presentation of analysis by distinguishing only between ``definite'' and ``indefinite'' quantifiers: the former govern domains for which a proof of consistency is available and secures the use of the law of excluded middle; the latter govern those for which there isn't, e.g.\ the real numbers. Lorenzen states in his foreword to \emph{Differential und Integral} (1965) that he is faithful to Weyl's approach of \emph{Das Kontinuum} in this simplification.

This history motivates a number of mathematical and philosophical issues about predicative mathematics: how does Weyl's interest into Lorenzen's operative mathematics fit with his turn to Brouwer's intuitionism as expressed in ``Über die neue Grundlagenkrise der Mathematik'' (\citeyear{weyl21})? Why does Lorenzen turn away from his language levels and how does this turn relate to Weyl's conception of predicative mathematics? What do Lorenzen's conceptions of mathematics reveal about Weyl's conceptions?
\[*\ *\ *\]

We propose a timeline for this history and describe briefly the issues related to it. Unless stated otherwise, the translations are ours. 

\section{1909--1918. \emph{Das Kontinuum}}
\label{sec:das-kontinuum}

The first three items of our timeline emphasise main themes of Weyl's monograph: real numbers and the problem of quantifying over them, the precedence of properties over sets.

\paragraph{1909. The rôle of Poincaré.}

According to Leon \citet[page~55]{chwistek35}, Poincaré's 1909 conferences in Göttingen \citep[see][]{poincare10} played a seminal rôle for Weyl's conception of real numbers in \emph{Das Kontinuum}.
\begin{quotation}
[Henri Poincaré] very positively asserted that real numbers do not form
a determined class and that it is meaningless to speak of all
real numbers. In the spring of 1909 I heard him deliver a
lecture at Göttingen, in which he expressed his point of view
very emphatically although in bad German. The youthful and
brilliant Hermann Weyl, an acknowledged disciple of Hilbert,
was present at this lecture. Undoubtedly it influenced his
conception of real numbers. \cite[pages~78-79.]{chwistek48}
\end{quotation}

\paragraph{1918. Criticism of sets as basic objects of mathematics.} We cite from \citealt[Chapter~1]{weyl18}, in the translation by \citet{pollardbole87}.

\begin{quotation}
  No one can describe an infinite set other than by indicating
properties\label{pred} which are characteristic of the elements of the set.
And no one can establish a correspondence among infinitely
many things without indicating a rule, i.e., a relation, which
connects the corresponding objects with one another. The
notion that an infinite set is a ``gathering'' brought together by
infinitely many individual arbitrary acts of selection, assembled
and then surveyed as a whole by consciousness, is nonsensical. \citep[page~23.]{weyl87}.
\end{quotation}

\citet{pollard05} investigates how this criticism stems from Weyl's reading of Fichte.

\paragraph{1918. Artificiality of types in analysis.}

\begin{quotation}
[\dots] in the iteration of the mathematical
process \emph{the two principles of closure~\(5\) [``filling in''] and~\(6\) [``there is\dots''] are to be applied only to
blanks which are affiliated with a basic category}.\textsuperscript{24}

\dotfill\smallskip

A ``hierarchical'' version of analysis is artificial and useless.
It loses sight of its proper object, i.e., number (cf.\ note~24).
Clearly, we must take the other path---that is, we must restrict
the existence concept to the basic categories (here, the natural
and rational numbers) and must not apply it in connection with
the system of properties and relations (or the sets, real numbers,
and so on, corresponding to them). In other words, the only
natural strategy is \emph{to abide by the narrower iteration procedure}.
Further, only this procedure guarantees too that all concepts and
results, quantities and operations of such a ``precision analysis''
are to be grasped as idealizations of analogues in a mathematics
of approximation operating with ``round numbers.'' This is of
crucial significance with regard to \emph{applications}. So a proposition
such as the one mentioned above, that every bounded set of real
numbers has a least upper bound, must certainly be abandoned.
But such sacrifices should keep the path ahead clear of
confusion.\textsuperscript{29}

\dotfill\medskip

\noindent\hfil
  \textbf{Notes to Chapter~1}

\dotfill\smallskip

24. The primary significance of the narrower procedure is most clearly
conveyed by the following observation: The objects of the basic categories
remain uninterruptedly the genuine objects of our investigation only when
we comply with the narrower procedure; otherwise, the profusion of
derived properties and relations becomes just as much an object of our
thought as the realm of those primitive objects. In order to reach a decision
about ``finitary'' judgments, i.e., those which are formed under the
restrictions of the narrower procedure, we need only survey these basic
objects; ``transfinitary'' judgments require that one also survey all derived
properties and relations.

\dotfill\smallskip

29. In science there are no ``commandments'', just ``laws''. So here too,
using the term ``there is'' in connection with objects which do not belong
to the basic categories should certainly not be ``forbidden''. It is, of course,
entirely possible (and permissible) to adopt the broader procedure; but if
this is done, let it be done in a non-circular way. 
\citep[pages~30, 32, 120; we have slightly amended the translation because we do not endorse the choice made in its note~37.]{weyl87}\end{quotation}








Note that the reedition \citealt{weyl60} has been reviewed by \cite{lorenzen60}. \citet{lorenzen86b}, in one of his last articles, starts with Weyl's \emph{Continuum} for a historical and philosophical discussion of the possibility of theorising the continuum.

\section{1947--1951. Lorenzen discovers Weyl as a predecessor}

The next four items of our timeline describe how Lorenzen's operative mathematics take shape and register Weyl as predecessor.

\paragraph{1947. First mention of Weyl's \emph{Continuum} to Lorenzen by Wilhelm Ackermann.}


Let us cite from a letter by Wilhelm Ackermann to Lorenzen dated 21 May 1947 (\foreignlanguage{german}{Paul-Lorenzen-Nachlass}, \foreignlanguage{german}{Philosophisches Archiv} of \foreignlanguage{german}{Universität Konstanz}, PL 1-1-117), whose subject is Lorenzen's manuscript ``Zur Neubegründung der Mathematik [On the new grounding of mathematics]'' (Hs.~974:155; also, with an introduction, in the Gottfried-Köthe-Nachlass at \foreignlanguage{german}{Universitätsarchiv Göttingen}, Cod.\ Ms.\ G.~Köthe M~10), in which he introduces his idea of language levels (see page~\pageref{levels}) as a ``calculus of regions''; we translate ``Kalkül'' by (logical) ``calculus'', whereas \citet[page~3]{lorenzen55} proposes the translation ``formal system''. 



\begin{quotation}
  [\dots] With regard first to the setup of your calculus of regions in general, which I can pretty much imagine according to your indications, I wish to encourage you to carry it out in detail as soon as possible. I also hold the view that a setup of mathematics free from contradiction (i.e.\ verifiably free from contradiction) must lie in the direction indicated by you or in a similar one, and the more so that I arrive at similar conceptions in my attempts to set up mathematics out of a type-free system of axioms that is verifiably free from contradiction. I do not believe any more that one may achieve a proof of freedom from contradiction for classical analysis.

  I am not as optimistic as you on the question whether the differences of such a setup with respect to classical analysis be merely irrelevant. It seems to me to the contrary that such a setup would not move that far away from that of Weyl (before he joined entirely the intuitionists). Consider first the theorem on the supremum of a bounded set of real numbers. It is clear to me that the theorem may be proved in the following form: every bounded set of real numbers over the \(n\)-th region admits a real number of the \(n+1\)-st region as supremum. But we are no longer talking about a real number per se and the set of real numbers seems to dissolve into sets of real numbers over the different regions. The same holds for the mean value theorem of continuous functions. You say at this point: the mean value theorem can be brought into the classical form if we form the union of all regions and define a real function as a definite sequence~\(f\), for which \(f_n\)~is a real function of the \(n\)-th~region such that \(f_{n+1}\)~is an extension of~\(f_n\). What do you mean here by the union of all regions? Do you mean the union of all finite regions? Then the defined real function belongs in my opinion to the region~\(\omega\). Or should one understand that you wish the regions to be extendable indefinitely upwards, and then form the union of all regions, i.e.\ with arbitrary ordinal? In this case you would not move inside the domain of a calculus verifiably free from contradiction. For the calculus of regions should be grounded in the freedom from contradiction of classical arithmetic. But there is no proof of freedom from contradiction of classical arithmetic per se. Namely, one may always state purely number-theoretic deductions that are not captured any more by a proof of freedom from contradiction at hand.\footnote{\selectlanguage{german}Was nun zunächst den Auf\bs
    bau Ihres Regionenkalküls im allgemeinen anbetrifft, von dem
    ich nach dem, was Sie schreiben, mir so ziemlich eine Vorstel\bs
    lung machen kann, so möchte ich Ihnen möglichst zureden, die\bs
    sen nun möglichst bald im einzelnen auszuführen. Ich bin auch der
    Ansicht, dass ein widerspruchsfreier (d.\ h.\ nachweislich wider\bs
    spruchsfreier) Aufbau der Mathematik in der von Ihnen angegebenen
    oder in ähnlicher Richtung liegen muss, und bin das um so mehr, 
    als ich bei meinen Versuchen, aus einem nachweislich wider\bs
    spruchsfreien typenfreien logischen Axiomensystem die Mathe\bs
    matik aufzubauen, auf ähnliche Begriffsbildungen komme. Dass
    sich ein Widerspruchsfreiheitsbeweis für die klassische Analysis 
    erzielen lässt, daran glaube ich nicht mehr.

    Nicht ganz so optimistisch wie Sie bin ich in der
    Frage, ob die Unterschiede eines derartigen Aufbaus gegen\bs
    über der klassischen Analysis nur unerheblich seien. Im Gegen\bs
    teil scheint es mir, dass ein derartiger Aufbau sich doch nicht
    so weit von dem von Weyl (bevor er sich ganz den Intuitionisten 
    anschloss) entfernt. Da wäre zunächst der Satz über die obere
    Grenze einer beschränkten Menge von reellen Zahlen. Es ist mir
    ohne weiteres klar, dass sich der Satz in der folgenden Form 
    beweisen lässt: Zu jeder beschränkten Menge von reellen Zahlen
    über der $n$-ten Region gibt es eine reelle Zahl der $n+1$-ten Region
    als obere Grenze. Aber von einer reellen Zahl schlechthin ist 
    dann doch keine Rede mehr, sondern die Menge von reellen Zah\bs
    len erscheint ein für allemal aufgelöst in Mengen von reellen
    Zahlen über den verschiedenen Regionen. Das gleiche gilt für
    den Zwischenwertsatz der stetigen Funktionen. Sie sagen an der
    Stelle: Der Zwischenwertsatz kann auch auf die klassische Form 
    gebracht werden, wenn wir die Vereinigung aller Regionen bilden
    und als reelle Funktion eine definite Folge~$f$ definieren, für
    die $f_n$~eine reelle Funktion der $n$-ten Region ist, sodass $f_{n+1}$ 
    eine Fortsetzung von~$f_n$ ist. Was meinen Sie hierbei mit der 
    Vereinigung aller Regionen? Meinen Sie die Vereinigung aller 
    endlichen Regionen? Dann gehört meiner Ansicht nach die de\bs
    finierte reelle Funktion der Region~$\omega$ an. Oder ist das so
    zu verstehen, dass Sie die Regionen nach oben unbegrenzt er\bs
    weiterungsfähig lassen möchten, und nun die Vereinigung aller 
    Regionen, d.\ h.\ solcher mit beliebiger Ordnungszahl bilden wol\bs
    len. In diesem Falle würden Sie sich nicht mehr in dem Bereiche 
    eines nachweislich widerspruchsfreien Kalküls bewegen. Denn
    der Regionenkalkül soll sich ja auf die Widerspruchsfreiheit
    der klassischen Arithmetik gründen. Es gibt aber nun keinen
    Widerspruchsfreiheitsbeweis für die klassische Arithmetik 
    schlechthin\barre{, ebenso wenig wie es ein Axiomensystem für die
      klassische Arithmetik schlechthin gibt}. Es lassen sich näm\bs
    lich immer rein zahlentheoretische Schlüsse angeben, die durch
    einen vorliegenden Widerspruchsfreiheitsbeweis nicht mehr er\bs
    fasst werden.}
\end{quotation}

\paragraph{1947. Second mention of Weyl's \emph{Continuum} to Lorenzen by Paul Bernays.}

Let us cite from a letter by Paul Bernays to Lorenzen dated  1 September 1947 (PL 1-1-112; a carbon copy lies in ETH-Bibliothek, Hochschularchiv, Hs.~975:2955).

\begin{quotation}
  Concerning the narrower calculus of types, it certainly has something principled to itself. But for the actual mathematical usage, it is yet quite unwieldy; and one will generally try to get along without a calculus of types. A formalism roughly equivalent to the narrower (ramified) type calculus was put up by Hermann Weyl in his writing ``Das Kontinuum''. I would think that the method of your proof of freedom from contradiction might be carried over to a formalism of about this kind.\footnote{\selectlanguage{german}Was den engeren Stufenkalkul anbelangt, so hat dieser ge\bs
wiss manches Prinzipielle für sich. Für den faktischen mathemati\bs
schen Gebrauch ist er aber doch wohl recht schwerfällig; und man
wird wohl überhaupt versuchen, ohne einen Stufenkalkul auszukom\bs
men. Ein dem engeren (verzweigten) Stufenkalkul ungefähr äquiva\bs
lenter Formalismus ist ja von Hermann Weyl in seiner Schrift
``Das Kontinuum'' aufgestellt worden. Ich möchte denken, dass auf ei\bs
nen Formalismus etwa dieser Art die Methode Ihres Wf.-Beweises
sich übertragen lassen sollte.}
\end{quotation}









\paragraph{1948. First inscription by Lorenzen of himself into the tradition of Weyl's \emph{Continuum}.}


In his letter to Bernays dated 15 August 1948 (Hs.~975:2961), Lorenzen writes about his enclosed manuscript ``Ein `trivialer' Widerspruchsfreiheitsbeweis für die Arithmetik [A `trivial' proof of freedom from contradiction for arithmetic]''  (Hs.~974:148) and concludes as follows.

\begin{quotation}
  I hope that in this form the question of freedom from contradiction, that in my opinion has disquieted mathematics long enough, has found a satisfactory answer. Here one has now firm ground for building up a constructive analysis (e.g.\ in the sense of Weyl), and then one will be able to examine also the axiomatic analysis as of its usability.\footnote{\selectlanguage{german}Ich hoffe, dass in dieser Form die Frage der Wf., die
  m.\ E.\ lange genug die Mathematik beunruhigt hat, eine befrie\bs
  digende Lösung gefunden hat. Hier hat man jetzt festen Boden,
  um eine konstruktive Analysis (z.\ B.\ im Sinne von Weyl) auf\bs
  zubauen, und dann wird man in aller Ruhe auch die axiomatische
  Analysis auf ihre Brauchbarkeit untersuchen können.}



\end{quotation}

\paragraph{1950. Lorenzen's relaunch of Weyl's approach.}

In his article ``Konstruktive Begründung der Mathematik'', \citet{lorenzen50} formulates his relaunch of Weyl's approach of 1918. We translate ``Begründung'' by ``justification'', but the word equally means ``foundation'', and ``establishment'' might be a good compromise; John Bacon uses the translation ``grounding'' in \citealt{MR0346097} (see the citation on page~\pageref{grounding}).

\begin{quotation}
  As most suggestive direction for further work results herefrom: first to develop in detail what has been summarised here under the name ``constructive set theory''. To the contrary of the intuitionistic attempts, one may use now throughout---after justifying logic unobjectionably---the excluded middle. The so uncomfortable restriction to ``decidable properties'', ``computable real numbers'', etc.\ is not necessary anymore. An approach in the direction intended here has been made by H. Weyl already in 1918---it had to fail because the justification of the excluded middle was still missing at that time.

  It is of course understood that the consolidation of a constructive analysis is of interest also independently of the goal of proving the usability of the traditional axiomatisations. New possibilities open up, e.g.\ for computing formally with infinite processes, independently of the concept of real number and convergence.

  The exploitation of such possibilities as the transfinite iteration for the process of set formation yields perhaps enough comfortable calculi to make the contemporary axiomatisations simply dispensable. (Pages~165--166.)
\end{quotation}

The very day he submits this article, 13 February 1950, he writes the following in a letter to Ackermann.
\begin{quotation}
  It seems thus more sensible to investigate the possibility of constructive set theories mindless of ``axioms''---instead of making random guesses among possible axioms with so-called plausibility considerations. The logic of ramified types (without reducibility) and Weyl's analysis are after all just the first attempts---if the distinction of types is dropped and only the ``ramification'' is retained (I call this then ``levels'')---then one obtains a further possibility which, as far as I see, leads to an analysis that practically does not differ from the non-formalised so-called naive analysis, as it is for instance common in lectures today. \citep[page~196.]{ackermann83}
\end{quotation}

In the follow-up article, submitted five months later, ``\german{Die Widerspruchsfreiheit der klassischen Analysis}'', Lorenzen emphasises as follows.

\begin{quotation}
  The problem of constructing real numbers has to my knowledge first been
  solved by H. Weyl (\emph{Das Kontinuum}, 1918). Admittedly, Weyl uses the arithmetic and logical rules without justification, in particular also the excluded middle. As a result of Brouwer's criticism, the intuitionistic analysis has then emerged, in which the fundamental theorems of classical analysis do not hold any more by the ban of the excluded middle and the restriction to a very narrow concept of function.

  Yet, after a constructive justification\textsuperscript{1} of arithmetic and logic (incl.\ excluded middle), a construction of real numbers may also be carried out, for which all fundamental theorems of classical analysis hold. The freedom from contradiction of classical analysis is contained therein. \citep[page~1.]{MR0043748}

  \vspace*{2.6pt}\noindent\hbox{\vbox{\footnoterule}}\vspace*{-5.6pt}

  \noindent{\footnotesize\rule{0pt}{\footnotesep}\textsuperscript{1} Cf. \citealt{lorenzen50}.}

\end{quotation}

He proceeds by explaining the process of levels\label{levels} of language and concludes.
\begin{quotation}
  Yet one does not obtain the classical analysis with the language over the rational numbers, for---if one defines the real numbers by an equality relation between sets of rational numbers (cf.\ §~2)---then no sets and functions of real numbers (but only those of rational numbers) are defined. To this end the construction of a language must be repeated: first the language over the previous statements, sets, functions, and real numbers (individuals of the 1st~level), then again a language over the thus obtained individuals of the 2nd~level, etc. (These levels must not be mistaken for the types of the naive and axiomatic set theory.)

  We call the union of all levels of finite height the 0th hyperlevel. The construction of a language over this one yields further hyperlevels whose individuals we call hyperstatements, hypersets, hyperfunctions, and hyperreal numbers. However, if one restricts to levels, then the classical completeness principle:
\[\text
  {``Each set of real numbers admits a real number as infimum''}
\]
becomes provable, although there are hypersets whose infimum is a hyperreal number. In §§~3--6 we obtain by an appropriate exclusion of hyperlevels all fundamental theorems of classical analysis (cf.\ for instance \german{Haupt--Aumann, \emph{Differential- und Integralrechnung~I}}). It\label{lines} is thus not the constructive analysis which is ``too narrow'' with regard to classical analysis: the closedness featured by classical analysis arises just by a restriction of the construction means. \citep[page~3.]{MR0043748}

\end{quotation}

\section{1955. Reception of Lorenzen's work by Weyl}

The next two items of our timeline tell how Weyl gets acquainted with Lorenzen.

In \citeyear{lorenzen55}, \citeauthor{lorenzen55} publishes his book \emph{Einführung in die operative Logik und Mathematik} that contains a comprehensive presentation of operational mathematics as a new foundation for mathematics. The term ``operational'' has been devised to indicate that the essential point is ``the operating according to rules'' \citep[page~4]{lorenzen55}; the alternative translation ``operative'' supersedes it only at a later point.

The reviews of the book are mostly positive or even enthusiastic (\citealp{kaminski55,robinson56,heinemann56,aubert56,becker57,frey57,skolem57,weizsaecker57}, \citealp[Chapter~III, §~7]{fraenkelbarhillel58}, \citealp{meigne70,heitsch71,demuth72,gardies72}); \citet{craig57} notes ``major though probably corrigible errors''; \citet{stegmueller58} wishes for ``a clear demarcation of the principles that are declared admissible for the reasonings in terms of content (the metatheoretical reasonings)'' and ``an elucidation of the strange, frequently encountered interlacing of object-linguistic and metalinguistic operations'';  the very detailed criticism and scepticism of \citet{mueller57} in the \emph{Zentralblatt} is puzzling because the preface of \citealt{lorenzen55} acknowledges his ``many good advices at the drafting of the manuscript and at the revisions'' (compare also \citealp{mueller57c}); four years later, this journal publishes another review by \citet{gericke61} that vindicates Lorenzen's book: ``an act of will is required not to measure the author's lines of thought by standards that are derived from other foundations. This difficulty can easily become the source of misconceptions''.

The theory proposed in \citealt{wang54} is ``essentially equivalent to the body of methods which Lorenzen
admits and applies'' \citep[page~78]{wang55}.

\paragraph{1955. ``Nachtrag Juni 1955'' to \citealt{weyl21}.}
\label{sec:nachtrag-}

In June 1955, \citeauthor{weyl56} writes an addendum to the inclusion of ``Über die neue Grundlagenkrise der Mathematik'' into his \emph{Selecta} (\citeyear{weyl56b}), in which he presents \citealt{lorenzen55} as ``most viable way out of the difficulties'' raised by the foundation of mathematics. 
It is translated in \citealt[page~12]{heinzmann21}.

Weyl's lines have impressed in particular Ignacio Angelelli. He reminds them as follows in his reviews of translations of \emph{Das Kontinuum} and of an article about it.
\begin{quotation}
  Readers who are interested in the philosophy of mathematics should take into account Weyl's later remarks on Paul Lorenzen's work as representing the best furtherance of his program (gangbarster Weg) \citep{angelelli91}.

  One misses any reference [\dots] to Paul Lorenzen, whose work was celebrated by Weyl in the mid 1950's as ``the most feasible approach to his program'' \citep{angelelli97,angelelli98}.
\end{quotation}
Note that Angelelli, a student of Józef Bocheński, spent the year 1965-1966 at the university of Erlangen for a Humboldt postdoctoral fellowship under the direction of Lorenzen \citep[see][]{legris20}.

\paragraph{1955. Weyl's letter to Lorenzen.}\label{sec:la-lettre-de}

The whole known correspondence between Lorenzen and Weyl lies in the \foreignlanguage{german}{Paul-Lorenzen-Nachlass} (PL 1-1-2).

Transcriptions of Weyl's letter dated 23 September 1955 and translated below lie in the  Helmut-Hasse-Nachlass (Cod.\ Ms.\ H. Hasse 1:1529, Beil.~16), in the Gottfried-Köthe-Nachlass (Cod.\ Ms.\ G. Köthe A~358, Beil.), and in the Josef-König-Nachlass (Cod.\ Ms.\ J. König 186) at \foreignlanguage{german}{Universitätsarchiv Göttingen} and show that Lorenzen communicated about it to his teachers: Hasse was his Ph.D.\ advisor; Köthe introduced him to lattice theory and was ready to replace Wolfgang Krull as his habilitation advisor \citep[see][page~241, note~50]{neuwirth21c}; König was his teacher of philosophy in Göttingen.

This letter has been presented before in \citealt{thiel00} and subsequently (in excerpts) in \citealt{schlaudt14}.





\begin{quotation}
  \noindent Hermann Weyl, Zürich\hfill 23 September 1955\nopagebreak

  \noindent\hfil Dear Mr. Lorenzen,\nopagebreak\medskip

  Your letter dated 17.9.\ must have crossed mine dated 18.9. In the meantime I have concerned myself in depth with your book and \emph{am most profoundly impressed by it}. Your ``operative'' standpoint actually seems to combine in the most natural and fortunate way formal construction inside the calculus with considerations in terms of content that lead to the insight of derivability, admissibility, etc. It is actually very simple, but yet surprising, how the possibility of inversion thus turns up! And you are able to circumnavigate the cliff of the excluded middle by means of Kolmogorov's and Heyting's thought of reinterpreting \(p \lor q\) into \(\dng{p \lor q}\). I naturally cannot peruse a book as yours at one go. But whenever I was returning to the book after cogitating your treatment of my own, I was enchanted how carefully and precisely you formulate everything. When I think in what restricted form I have afforded the construction of relations in my little book of 1918 (by the principle of substitution and iteration), I cannot help being surprised by the free way in which you handle the inductive definition of relations~(\(\mathsf{I}_\varrho\)). I have not yet coped with this, but neither do I see what one could object against this broad-mindedness. I have noticed all the same that you yourself have found it good with regard your first account in 1951 to restrict the induction schemes through the requirement ``founded'' and ``separated'' (for the sake of definiteness). I do not remember today anymore precisely why I shied away from the iteration of the ``mathematical process'' in 1918, for you it is essential to iterate, even up to a limit ordinal.

  In spite of the so much further stretched frame, I naturally recognise the methodical kinship of your treatment for the justification of analysis with my much more limited approach of 1918---and am pleased by it! At the time, my belief into the truth of the excluded middle in terms of content was still unimpaired in the domain of the natural numbers. Brouwer had to come to free us from it, and in this respect you are also liable to Brouwer. Ultimately also to Hilbert: he emphasised the formalism of the calculus, his program of a proof of freedom from contradiction was perspicuous and promising in its first stages, one had first to experience the difficulties that show that this attack is as good as hopeless before one was ready to go other ways---\emph{your} way.

  I feel that I see at last again a clear sky after long years of resignation. I am happy that I have lived to see this. Please accept once again my warmest congratulations for the work done by you!---You hand this book over to me as present for my (soon due) 70th birthday. It has already arrived this morning. I believe that one could not have made me a finer present. I wholeheartedly thank you for this!\strut

  \hfil\hfil\begin{tabular}{c}
    With best greetings and wishes\\
    your\\
    Hermann Weyl\footnotemark
  \end{tabular}
\footnotetext{\selectlanguage{german}\noindent Hermann Weyl, Zürich\hfill 23.\ September 1955\nopagebreak

    \noindent\hfil Lieber Herr Lorenzen,\nopagebreak\medskip

    Ihr Brief vom 17.9. muss sich mit dem meinigen vom 18.9.\ gekreuzt haben.
    In der Zwischenzeit habe ich mich eingehend mit Ihrem Buch beschäftigt und
    \emph{bin aufs tiefste davon beeindruckt}. Ihr ``operativer'' Standpunkt scheint wirk\bs
    lich auf die natürlichste und glücklichste Weise formale Konstruktion inner\bs
    halb des Kalküls mit inhaltlichen Überlegungen, die zur Einsicht in Ableit\bs
    barkeit, Zulässigkeit etc.\ führen zu kombinieren. Es ist ja wirklich sehr
    simpel, aber doch überraschend, wie sich so die Möglichkeit der Inversion
    herausstellt! Und die Klippe des tertium non datur können Sie mittels des
    Kolmogoroff-Heyting'schen Gedankens der Umdeutung von \(p \lor q\) in \(\dng{p \lor q}\) um\bs
    schiffen. Ich kann natürlich ein Buch wie das Ihre nicht auf einen Rutsch
    durchlesen. Aber immer, wenn ich nach eigenem Nachdenken über Ihr Verfahren
    zu dem Buch zurückkehrte, fand ich mich entzückt davon, wie sorgfältig und
    präzis Sie alles formulieren. Wenn ich daran denke, in wie eingeschränkter
    Form ich mir in meinem Büchlein von 1918 die Konstruktion von Relationen er\bs
    laubt hatte (durch das Substitutions- und Iterations-Prinzip), so kann ich
    nicht anders, als erstaunt sein über die freie Weise, in der Sie die induktive
    Definition von Relationen (\(\mathsf{I}_\varrho\)) handhaben. Das habe ich noch nicht ganz be\bs
    wältigt, aber ich sehe auch nicht, was man gegen diese Weitherzigkeit ein\bs
    wenden könnte! Immerhin bemerkte ich, dass Sie selbst es für gut fanden,
    gegenüber Ihrer ersten Darstellung 1951, die Induktionsschemata durch die For\bs
    derung ``fundiert'' und ``separiert'' einzuschränken (um der Definitheit willen).
    Ich erinnere mich heute nicht mehr genau, warum ich 1918 vor der Iteration
    des ``mathematischen Prozesses'' zurückschreckte, für Sie ist es wesentlich zu
    iterieren, sogar bis zu einer Limeszahl.

    Trotz des so viel weiter gespannten Rahmens erkenne ich natürlich die
    methodische Verwandtschaft Ihres Verfahrens zur Begründung der Analysis mit
    meinem viel limitierteren Ansatz von 1918 -- und bin erfreut darüber! Damals
    war bei mir noch der Glaube an die inhaltliche Wahrheit des tertium non datur
    im Gebiete der natürlichen Zahlen unerschüttert gewesen. Brouwer musste kom\bs
    men, uns davon zu befreien, und insofern sind Sie auch Brouwer verpflichtet.
    Schließlich auch Hilbert: er betonte den Formalismus des Kalküls, sein Pro\bs
    gramm eines Beweises der Widerspruchslosigkeit war einleuchtend und in den ersten
    Stadien versprechend, man musste erst die Erfahrung der Schwierigkeiten
    machen, die diesen Angriff so gut wie hoffnungslos erscheinen lassen, bevor
    man bereit war, andere Wege -- \emph{Ihren} Weg zu gehen.

    Mir geht es so, dass ich nach langen Jahren der Resignation endlich wie\bs
    der offenen Himmel sehe. Ich bin froh, dass ich das noch erlebt habe. Nehmen
    Sie noch einmal meine herzlichsten Glückwünsche entgegen für das von Ihnen
    vollbrachte Werk! -- Sie überreichen mir Ihr Buch als Geschenk zu meinem
    (bald fälligen) 70.~Geburtstag. Es ist heute morgen schon eingetroffen. Ich
    glaube, man hätte mir kein schöneres Geschenk machen können. Ich danke Ihnen
    von ganzem Herzen dafür!\strut

    \hfil\hfil\begin{tabular}{c}
      Mit den besten Grüßen und Wünschen\\
      Ihr\\
      Hermann Weyl
    \end{tabular}}
\end{quotation}
\section{1955-1958. Weyl's invitation of Lorenzen to the Institute for Advanced Study}\label{sec:linv-de-lorenz}

The next six items of our timeline document the missed opportunity of our two protagonists meeting and Gödel's irritation by Lorenzen's operative mathematics.

This invitation is documented in the file ``Lorenzen, Paul, 1955--1958'' of the School of Mathematics records, 00110, of the Institute for Advanced Study in the Shelby White and Leon Levy Archives Center.




\paragraph{1955. Letter from Weyl to Selberg.} This letter is not dated.

\begin{quotation}
  Dear Atle, I am writing to you to suggest as a candidate for
  the Institute, not for this, but maybe for the next academic year,
  Prof.\ \emph{Paul Lorenzen} from the University Bonn (address: Carl-Schurz\Bs
  Colleg, Kaiserstr.~57, Bonn). He is, of course, German, 40~years
  old, married and has one daughter, age~13. He is best known for
  his research in mathematical logic and foundations of mathematics.
  But he has also published a number of good papers in algebra.

  Quite recently a book of his ``Einführung in die operative
  Logik und Mathematik'' appeared in the Gelbe Sammlung (Sprin\-ger-Ver\bs
  lag). Of this book both van der Waerden and I think very highly.
  I am reading it at the moment with considerable enthusiasm. It seems
  to me that, after Brouwer's intuitionism and Hilbert's formalism
  and Gödel's debacle, this may actually show the right and best way
  out of the dilemma, and give us---if not the consistency proof of
  one of the formal axiomatic systems---at least all the essential
  theorems of analysis. On Lorenzen's ``operative'' standpoint Gödel's
  discovery loses completely its disquieting character.

  Whether my judgment is right or wrong: Lorenzen is certainly
  an investigator in his field of high originality and one whom it is
  worthwhile to bring over for a visit to the States. I do not know
  him personally, but most mathematicians here do, and they all tell
  me that he is a very lively and agreeable fellow.

  There is an additional reason why I could wish Lorenzen to
  have an opportunity to visit America: the younger Hirzebruch (and
  not he) got the vacant full professorship in Bonn, and H.\ will go
  there after his year at Princeton University. L.\ is merely what
  they call an apl.\ (meaning, I guess, ``außerplanmässiger'') professor
  in Bonn.

  \hfill H. Weyl
\end{quotation}
\paragraph{1955. School minutes of the School of Mathematics of the Institute for Advanced Study, 26 October.}
\label{sec:scho-minut-scho}

\begin{quotation}
  \begin{itemize}
  \item[3.] The case of \emph{Paul Lorenzen} of the University of Bonn was
    brought up by Professor Selberg. It is understood that Professor Gödel
    will examine Lorenzen's papers and report on them.
  \end{itemize}
\end{quotation}

\paragraph{1955. Gödel's first report on Lorenzen.}
The following draft by Gödel are taken from a file entitled ``Lorenzen, letter, drafts, and discussion notes: 1958 and undated'', Kurt Gödel Papers, Personal and Scientific Correspondence, Lorenzen, Paul, Box 2b, Folder 103, documents number 011467--011468, C0282, Manuscripts Division, Department of Special Collections, Princeton University Library. This draft is an elaboration of a sketch (011464, 011466, 011469). We have typeset our completion of Gödel's abbreviations in oblique type.

According to Jan von Plato, they contain ``perhaps the nastiest passages on scientific matters in the thousands of pages of the Kurt Gödel Papers''.

\begin{quotation}
First: also purely algebraic papers which I have not read bec\sub.{ause} I
don't feel comp\sub.{etent} to judge about them\add.
His papers in the found\sub.{ations} refer to 3 different topics.

1. He has extended Weyl's constr\sub.{uctive} approach to analysis by intro\bs
ducing a hierarchy of order\add{s} of real nu\sub.{mbers}
\& he got as far as defining Leb\sub.{esgue} integral \& proving its
main properties. Also he proved Urysohn's theorem
about the introd\sub.{uction} of a metric in a topological space\add.
This work is very good as far as it goes. But the
res\sub.{ults} are not very prof\sub.{ound}. It is a well-known fact
that large parts of analysis do not use impred\sub.{icative} procedure\add s.
In many cases one only has to analyse the proofs in order to find it out.
Moreover this approach has its limits \& therefore it is entirely
unjustified to say (as Lorenzen does) that thereby the cons\sub.{istency}
of analysis has been proved.
In my op\sub.{inion} this whole appr\sub.{oach} goes in the wrong
direction.\footnote{In the sketch for this draft, Gödel is more explicit: he writes that ``the value of this whole
  approach to analysis is rather questionable'' and refers to three sentences of \citet[page~23]{weyl18}, copied in Gabelsberger short:

  ``1. An analysis with formation of types is artificial and useless.
  
  ``2. Every cell of this mighty organism is permeated by the poison of contradiction.

  ``3. We must restrict the existence concept to the basic categories.''\footnotemark}\footnotetext{%
  \begin{tabular}[t]{@{}l@{}}
    1. Eine Analysis mit Stufenbildung ist künstlich und unbrauchbar.\\
    2. Jede Zelle des gewaltigen Organ\add{ismus} der Analysis ist von diesem Gift durchsetzt.\\
    3. Der Existenzbegriff soll nur hinsichtlich der Grundkategorien angewendet werden.
  \end{tabular}}

2. The second topic of Lorenzen's work is a new method
for consistency proofs.\footnote{\label{lor1951}\citealt{lorenzen51wt}, see \citealt{coquandneuwirth23}.} This method constitutes a \uline{very}
import\sub.{ant} progress in this field. However this method is
really due to Novikoff who published it 7 years earlier
(1944) in English and a sketch of it even came out 10 years
earlier.\footnote{\citealt{novikoff39,novikoff43}.} Novikoff's proof is not quite satisfactorily
refereed
but still I would say by far the greater part of the credit
has to go to Novik\add{off}. Moreover Lorenzen also does
not explicitly give the \uline{really interesting} result obtainable by
this method and it is not trivial to extricate it from his proof. It was
Schütte who did that.

3. The third topic of Lorenz\sub.{en's} work is a foundation of arithmetic
on the basis of his operational
viewpoint. It cons\sub.{ists} in a mixture of form\add{alism} \& intuitionism.
It is vague\add: it does not dist\sub.{inguish}
between formulas as marks on paper \& as expr\sub.{essions} having meaning.
The proofs are not carried out at all
but nevertheless he claims to have given a satisfactory
found\add{ation} of arithmetic. It is true that it contains ideas which probably
could be carried out
\& possibly might lead to something interesting
but then this found\add{ation} of arithm\add{etic} would
change its character \& would rather lead to
a disproof of the operational viewpoint namely
it would show that if
one wants to give a satisf\sub.{actory} found\add{ation} for number th\sub.{eory}
(including a cons\sub.{istency} proof) it is not possible to confine oneself to cons\sub.{istently}
referring to the handling of symbols but has to use cert\add{ain} abstr\sub.{act}
conc\sub.{epts} such as function, implic\add{ation}, etc.
%

I am definitely opposed to inviting Lorenzen\add.
My judgement is based on his work on the found\add{ations.} He wrote
some purely alg\sub.{ebraic} pap\sub.{ers}
about ord\add{ered} gr\sub.{oups} \& lattice gr\sub.{oups}
which I have not read but they
are not the reason why Weyl wants to invite him\add.
As to the foundations there is I.~his absurd phil\sub.{osophy} which
also reflects very unfavorably
on his math\sub.{ematical} work namely 1.~It induces him to
make openly false statements about his math\sub.{ematical} work
1.~that he has used no transf\sub.{inite} ind\add{uction} where he evidently has used
it\footnote{Compare the account of accessibility in \citealt[note~c on page~272,]{goedel90} and \citealt[§~6.3]{coquandneuwirth23}} or to say he has proved the cons\sub.{istency} of anal\sub.{ysis} where
there is nothing of the kind. One even gets the impr\sub.{ession} that
these ass\sub.{ertions} are consciously wrong\add. 2.~Moreover which is even more serious
his phil\add{osophical} prejudices prevent him from doing the math\add{ematical} work properly. He takes the pos\add{ition}
that there is
nothing like math\sub.{ematical} cogn\sub.{ition} or a math\add{ematical} intuition
or evident math\add{ematical} axioms
but that everything in math\add{ematics} is due to conv\sub.{ention} \& to expediency.
This prejudice prevents him from
analysing what he uses in his metamath\sub.{ematical} proofs which
is of preeminent[?]\ imp\sub.{ortance}
in this field. Moreover
this has the cons\sub.{equence} that he is imprecise in his metamath\sub.{ematical} reasoning\add.

II.~Considering his papers one by one the best no doubt is the one in \add{The Journal of Symbolic Logic.\footnotemark[\getrefnumber{lor1951}]}
\end{quotation}

Note that in a letter dated 20 March 1956, Gödel asks von Neumann what he thinks about Lorenzen's ``efforts to base analysis on ramified type theory [\dots] up to the theory of the Lebesgue measure'' \citep*[pages~376--377]{goedel03b}.

\paragraph{1955. Research intended by Lorenzen in Princeton.}
\label{sec:rese-intend-lorenz}

\begin{quotation}
  \indent\indent Intended research.\bigskip

  I intend to continue the investigations I have started
  in my book. Above all I hope to achieve a progress in one
  of the following points:
  \begin{enumerate}[(1)]
  \item application of the operational method to functional
    analysis. In some lectures I have already developed
    an approach to additive set-functions and linear
    functionals.
  \item examination, as to how far the theorems of abstract
    topology may be reached ``operationally''.
  \item application of the operational method to metamathe\bs
    matics in the sense of Tarski and Robinson.
  \end{enumerate}

  \hfill
  \begin{tabular}{r}
    Paul Lorenzen\\
    Nov.~27, 55
  \end{tabular}
\end{quotation}


\paragraph{1956. School minutes of the School of Mathematics of the Institute for Advanced Study, 11 April.}

\begin{quotation}
  \begin{itemize}
  \item[5.] It was voted to offer \emph{Paul Lorenzen} a grant of \$4,700 and membership for the academic year 1957--58. [Professor Selberg will write informally first.]
  \end{itemize}
\end{quotation}

\paragraph{1958. Gödel's second report on Lorenzen.}
\label{sec:godels-report-lorenz}

A draft of the report below is in the Kurt Gödel Papers (011463).\bigskip

{\parindent0em\noindent
  To be returned to the Committee at the conclusion of the scholar's stay at his host in\bs
  stitution.\smallskip

  \hfil\makebox[0pt][c]{REPORT TO THE COMMITTEE ON INTERNATIONAL EXCHANGE OF PERSONS}\smallskip

  \hfil Conference Board of Associated Research Councils

  \hfil 2101 Constitution Avenue, Washington 25, D. C.\rlap{\hspace{4em}1957-8}\smallskip

  \leavevmode\rlap{Name of Visiting Scholar: \textb{Paul LORENZEN}}\hphantom{Sponsoring Institution: \textb{Institute for Advanced Study}\qquad}Country \textb{Germany}\smallskip
  
  Sponsoring Institution: \textb{Institute for Advanced Study}\qquad Field \textb{Mathematics}\smallskip
  
  Member of faculty most familiar with his work: \textb{Professor Kurt Gödel}\smallskip

  \leavevmode\rlap{Arrival date at the institution:}\hphantom{Arrival date in U.S.\ (if known): }\textb{9/12/57}\enspace Departure date from institution: \textb{3/30/58}\smallskip
  
  Arrival date in U.S.\ (if known): \textb{9/12/57}\enspace Departure date from U.S.\ (if known): \textb{3/30/58}\smallskip
  
  Forwarding address if scholar is still in the United States:\smallskip

  Academic Status of Visitor: \textb{Member of the School of Mathematics}\smallskip

  Has the scholar been primarily engaged in research \hbox to 3em{\hfil\textb{Yes}\hfil}; in teaching \hbox to 3em{\hfil\textb{No}\hfil}?\smallskip

  If teaching, the following information would be helpful to the Committee:\smallskip

  {\leftskip 1em Approximate number of students taught:\qquad Number of teaching hours per week:\smallskip

  Did he teach primarily graduate or undergraduate students?\smallskip
  
  Did he teach courses regularly offered in the curriculum?\qquad (On reverse side,
  please list titles of courses taught)\smallskip

  }If he undertook research please comment briefly on the nature of his research, indicat\bs
  ing its value (a)~to your institution or (b) to his own professional development.\smallskip

  \textb{He obtained a considerable improvement and strengthening of one of the axioms of his
    ``operationistic'' mathematics and with its help, gave a proof of the Cantor-Bendixson
    Theorem which is constructivistic in a generalized sense. He also conceived interesting
    ideas of a game-theoretical interpretation of logic.}\smallskip

  Did the visitor have any difficulties with the language? \qquad. If so, was this a
serious handicap to him in pursuing his professional program? \qquad; in his social
contacts? \hspace*{2em}. (A reply to this question is requested if this is a report on a
visitor from Asia.)\smallskip

Additional remarks: (Please use reverse side if needed.)\bigskip

  \noindent
  \begin{tabular}[b]{l}
    Date: \textb{June 25, 1958}\\[\smallskipamount]
    \begin{tabular}[c]{@{}l}
      EPL:ej\\
      Oct.~56
    \end{tabular}\qquad August 14, 1958
  \end{tabular}
  \hfill
  \begin{tabular}[b]{l}
    Signature: /s/ Kurt Gödel\\[\smallskipamount]
    Signature: \textb{Kurt Gödel}\\[\smallskipamount]
    Title: \textb{Director}
  \end{tabular}

  }

\section{1958--1965. Transition of Lorenzen to a simple distinction of definite and indefinite}
\label{sec:trans-simple-dist}

The last four items of our timeline describe how Lorenzen shifts to dialogues and indefinite quantifiers.

\paragraph{1958. The transition to dialogical logic.}

As Kuno \citet{lorenz21} tells, Lorenzen and Tarski met at the international symposium ``The axiomatic method with special reference to geometry and physics'' held at the university of
California at Berkeley from 26 December 1957 to 4 January 1958.
Lorenz witnesses that ``the further development of the operative to a \emph{dialogical} logic based on validity claims and their defense resp.\ rejection has not become a desideratum for Lorenzen before Tarski's doubts about the appropriateness of his usage of `definite'{}'' (private communication, 8 February 2022).

\citet{lorenz01} describes the problematic of the meaning of implication (called subjunction) in operative logic. It is based on \emph{admissibility}, so that one has to move one level up in the following sense:  if the rule leading from~\(A\) to~\(B\) is admissible in a calculus, i.e.\ if it does not extend the class of derivable statements, then one states that \(A\implies B\) at a metalevel. However, admissibility is undecidable. Lorenz writes: ``Especially subjunctions in the 
operative interpretation could not any more be called ‘definite’ as it had been 
the explicit intention of the operative approach'' (page~257). He proceeds with pointing out the incapacity of operative logic to capture the meaning of what it is to be a proposition and the solution provided by the concept of ``dialogue-definiteness''.

\citeauthor{schroederheister08} (\citeyear{schroederheister08}, §~3.2, \citeyear{schroederheister08b}, §~3.1) describes Lorenzen's theory of implication in his comparison of operative mathematics with proof-theoretic semantics.

\paragraph{1965. No language levels in Lorenzen's \emph{Differential und Integral}.}

The foreword to \citealt{MR0209410} underlines the simplification in abandoning the construction of language levels.

\begin{quotation}
    
    


    
  Since the appearance of my \emph{\german{Einführung in die operative Logik und
Mathematik}} (\citeyear{lorenzen55}), I have given a simplified account of the
grounding\label{grounding} of logic in \emph{Metamathematik} (\citeyear{lorenzen62}). I now put forth a
grounding of analysis which is likewise considerably simplified:
as was the case in H. Weyl's \emph{Das Kontinuum} (\citeyear{weyl18}), no construc\bs
tion of ``higher'' language levels is undertaken. For this reason
the book is dedicated to the memory of Hermann Weyl. \citep[page~\textsc{ix}.]{MR0346097}
\end{quotation}


\citet{feferman00} reports as follows on \citealt{MR0346097}: ``while significant portions of that are based on predicative grounds, it is not restricted to such''. Feferman does not justify his judgment; we presume that it targets the consideration of inductive definitions: as Laura \cite{crosilla22} points out, he is sticking to the ``classical approach'' of a ``predicativity given the natural numbers'' \citep{feferman05}.


\paragraph{1965. Bernays reviews Lorenzen's \emph{Differential und Integral}.}

At the request of \german{Akademische Verlagsgesellschaft}, the publisher of \emph{Differential und Integral}, Bernays writes a review (Hs~975:5329) on 5 November 1965, and in particular the following.
\begin{quotation}
  The specifically remarkable about this book from the foundational point of view consists in avoiding consistently the ``impredicative'' methods. This methodical restriction, demanded first at the beginning of our century by several French mathematicians, and for which on the one hand  Russell und Whitehead with their ramified theory of types, on the other hand Hermann Weyl in his writing ``Das Kontinuum'' provided a formal framework, is carried out here in such a way that the explicit put-up of a formal framework is not required, that it is indeed enough to tighten the usual methods in the sense of a \emph{constructive version of the concept of existence for real numbers}, so that the variety of real numbers needn't be considered as a determined delimited one, but can rather be treated as an ``indefinite'' totality. The careful execution of this program confirms the view harboured by theorists of foundations that in doing so no essential perceptible restriction of the methods and results of analysis occurs.\footnote{\selectlanguage{german}
%
%
%
%
%
  Das spezifisch Bemerkenswerte des Buches vom grundlagentheoreti\bs
  schen Standpunkt besteht in der konsequenten Vermeidung der "`imprädika\bs
  tiven"' Methoden. Diese methodische Beschränkung, wie sie zuerst im
  Anfang unseres Jahrhunderts von mehreren der französischen Mathematiker
  gefordert wurde und für welche einerseits Russell und Whitehead mit 
  ihrer verzweigten Stufen\-theorie, andrerseits Hermann Weyl in seiner
  Schrift "`das Continuum"' einen formalen Rahmen lieferten, wird hier in
  solcher Weise durchgeführt, dass es der expliciten Aufstellung eines
  formalen Rahmens für die Darstellung nicht bedarf, dass es vielmehr
  genügt, die üblichen Methoden im Sinne der \emph{konstruktiven Fassung des
    Existenzbegriffes für reelle Zahlen} zu verschärfen, sodass die Mannig\bs
  faltigkeit der reellen Zahlen nicht als eine bestimmt abgegrenzte in
  Anspruch genommen zu werden braucht, vielmehr als eine "`indefinite"'
  Gesamtheit behandelt werden kann. In der sorgfältigen Durchführung dieses
  Programmes bestätigt sich die schon von den Grund\-lagen\-theo\-re\-ti\-kern ge\bs
  hegte Ansicht, dass hierdurch die Methoden und Ergebnisse der Analysis
  keine wesentlich fühlbare Beschränkung erfahren.



  }
\end{quotation}

\paragraph{1969. Foreword to the second edition of Lorenzen' \emph{Einführung in die operative Logik und Mathematik}: indefinite quantifiers are more appropriate.}

In the foreword, \citet{MR0241272} apologises for keeping with an approach superseded by his later work.

\begin{quotation}
  Although I---understandably---hold the approaches of later works to be ``more appropriate'', e.g.\ a logic of dialogues instead of a logic of calculi, the usage of indefinite quantifiers instead of an explicit construction of language levels, this new edition contains the content of the 1st edition unchanged.

  The reader can thus carry out himself a comparison with my later works (cf.\ list of references).
  

\end{quotation}

The later works alluded to are \citealt{lorenzen62,MR0209410,MR0290929}.

\section{Focus on a few issues}
\label{sec:focus-few-points}

\paragraph{No presupposition of a totality in the constructive method.}

Let us cite \citealt{hilbertbernays34} in the translation of the \cite{bernaysproject03}: the axiomatic method involves a further assumption with respect to the constructive method, viz.\ the assumption of a fixed system of things. 
\begin{quotation}
  Another factor coming along in axiomatics in the narrowest sense is the \emph{existential form}. It serves to distinguish
  \emph{the axiomatic method} from the \emph{constructive} or \emph{genetic} method of founding
  a theory.\thinspace\textsuperscript{1} Whereas in the constructive method the objects of a theory
  are introduced merely as a \emph{family} of things,\thinspace\textsuperscript{2} in an axiomatic theory one
  is concerned with a fixed system of things (or several such systems) which
  constitutes a previously \emph{delimited domain of subjects} for all predicates from
  which the statement of the theory are constituted.
  
  Except in the trivial cases in which a theory has to do just with a finite,
  fixed totality of things, the presupposition of such a totality, of a ``domain
  of individuals'', involves an idealizing assumption joining the assumptions
  formulated in the axioms. \citep[pages~1--2.]{hilbertbernays34}

  \footnotesize\vspace*{2.6pt}\noindent\hbox{\vbox{\footnoterule}}\vspace*{-3pt}

  \noindent\rule{0pt}{\footnotesep}\textsuperscript{1} See for this comparison appendix~VI to Hilbert's Grundlagen der Geometrie: Über
    den Zahlbegriff, 1900.

    \noindent\textsuperscript{2} Brouwer and his school use the word ``species'' in this sense.
\end{quotation}

Compare the last lines of the quotation from \citealt{MR0043748} translated on page~\pageref{lines}.

This ``genetic'' aspect of predicative mathematics is described as ``logical reflection'' in \citealt{MR0104569}.
\begin{quotation}
  But functions and relations are not the objects of arithmetic.
They are the concepts used in speaking about numbers as the proper objects.
Now the transition from arithmetic to analysis is achieved by taking as
the objects of a new theory just these concepts of the old theory. Psycholo\bs
gically expressed, the focus of attention has to pass from the old objects,
the numbers, to the functions and relations as new objects. Let me call this
transition a ``logical reflection'', because the reflection to be performed is
on the concepts occurring in the theorems of the old theory. Or more briefly,
the object of the reflection is the language used so far. In view of this one
may be justified in calling it a ``logical'' reflection. This logical reflection,
of course, does not simply mean that we become aware of the concepts
used so far; it also means that we ask ourselves the question: which concepts
could ``possibly'' be used. Put into words without safeguarding against
future difficulties, one might ask for ``the class of all possible arithmetical
functions and relations''. With this question, with this logical reflection,
the transition from arithmetic to analysis begins. It is the same process
which Hermann Weyl in his book \emph{Das Kontinuum} (Leipzig, 1918) simply
called ``the mathematical process'' because of its fundamental importance
to modern mathematics. Greek mathematics differs from modern mathe\bs
matics just by not having achieved this logical reflection. \citep[page~244.]{MR0104569}.
\end{quotation}

\citet[§~4]{kahleoitavem21} discuss the article \citealt{MR0043748} and judge that because of its account of the ``mathematical process'' by language levels, ``one can hardly call the system under consideration `classical analysis'{}'' and that ``it also marks a significant turn of Lorenzen away from Hilbert’s philosophical basis''. These judgments seem to consider that such a philosophical basis must aim at justifying Dedekind's impredicative introduction of real numbers.


\paragraph{Definiteness and 
  inductive definitions.}
\label{sec:defin-induct-gener}

\citet{lorenzen55} introduces the concept of \emph{definite} proposition in order to ``specify the \emph{methods} with which the subject matter [of mathematics] may be investigated resp.\ recognised'' (page~4). In doing so, he leaves the notion of schematical operation undefined: ``To operate schematically with figures is familiar to everybody'' (page~9).
\begin{quotation}
  With a view of employing no unnecessary or arbitrary prohibition, the present essay leaves the methodical frame as large as possible. A limit that is insuperable for every part of mathematics that should  be considered as ``firm'', ``secure'' (or however one would like to name it), seems to me to lie in the propositions to be ``definite''.
  
  \dotfill\smallskip

  We give therefore the following inductive definition of ``definite'':
  \begin{enumerate}[(1)]
  \item Every proposition that is decidable by schematical operations is definite.
  \item If a definite concept of proof or of refutation is stipulated for a proposition, then the proposition itself is definite, more precisely proof-definite, resp.\ refutation definite.
  \end{enumerate}

  By the methodical demand of definiteness, the impredicative conceptions are of course excluded---as demanded already by Russell and Poin\-caré---, but quantifiers are not. If \(A(x)\) is definite, then let us stipulate as refutation for the proposition \(\bigland_xA(x)\) [for all~\(x\): \(A(x)\)]: a refutation of a proposition~\(A(x_0)\). For the proposition \(\biglor_xA(x)\) [for some~\(x\): \(A(x)\)], let us stipulate as proof: the proof of a proposition~\(A(x_0)\). \citep[pages~5--6.]{lorenzen55}
\end{quotation}

\citet[§~16]{lorenzen55} considers 
inductive definitions as providing operative counterparts to impredicative conceptions. He works out the example of the concept of convergent subsequence \(r_{l_*}\) of a sequence~\(r_*\) of rational numbers between~\(0\) and~\(1\) with the following proof.
\begin{quotation}
  We let \(l_1=1\). We halve then the interval~\(I_1=[0,1]\). If both halves contain infinitely many terms of the sequence~\(r_*\), then let~\(I_2\) be the left half, else the half that contains infinitely many terms. Let~\(l_2\) be the minimal number such that \(l>l_1\) and \(r_l\in I_2\). Subsequently \(I_2\)~is halved, \(I_3\)~is defined as one of the halves, then \(l_3\)~is determined analogously, etc. \citep[page~166.]{lorenzen55}
\end{quotation}
He proposes to defines the concept of an interval~\([a,a+2^{1-n}]\) containing infinitely many terms of~\(r_*\) as a relation \(n\Srel a\) by the following 
inductive rules:
\[
  \begin{aligned}
    &\begin{prooftree}
      \infer0{1\Srel 0}
    \end{prooftree}\\[3pt]
    &\begin{prooftree}
      \hypo{n\Srel a}\hypo{\bigland_{m_0}\biglor_{m>m_0}a\le r_m\le a+2^{-n}}\infer2{n+1\Srel a}
    \end{prooftree}\\[3pt]
    &\begin{prooftree}
      \hypo{n\Srel a}\hypo{\biglor_{m_0}\bigland_{m>m_0}(r_m<a\lor r_m> a+2^{-n})}\infer2{n+1\Srel a+2^{-n}}
    \end{prooftree}
  \end{aligned}
\]
As Lorenzen notes, these rules stipulate how a proposition~\(n\Srel a\) is to be proved, but it might well be that one does not obtain a proof of neither~\(2\Srel0\) nor~\(2\Srel2^{-1}\) (an instance of the law of excluded middle).

Then he defines the sequence of indices~\(l_*\) inductively by the following inductive rules.
\[
  \begin{aligned}
    &\begin{prooftree}
      \infer0{l_1=1}
    \end{prooftree}\\[3pt]
    &\begin{prooftree}\hypo{l_n=k}\hypo{n+1\Srel a}\hypo{m=\min\{\,l>k:a\le r_l\le a+2^{-n}\,\}}\infer3{l_{n+1}=m}
    \end{prooftree}    
  \end{aligned}
\]
Lorenzen concludes that the propositions \(l_n=k\) are definite, so that the concept of convergent subsequence is definite, regardless of its existence.

He introduces the following notation and terminology for the 
inductive definition of the relation~\(l_n=k\) between~\(n\) and~\(k\):
\[\left.
  \begin{aligned}
    &\,\to\,{1,1}\in\varrho\\
    {n,k}\in\varrho\,\land\,{n+1\Srel a}\,\land\, {m=\min\{\,l>k:a\le r_l\le a+2^{-n}\,\}}&\,\to\, {{n+1,m}\in\varrho}
  \end{aligned}\right\}=\mathfrak A(\varrho)\]
is the corresponding \emph{induction scheme} and the sought-after relation is
\(\mathsf{I}_\varrho\mathfrak A(\varrho)\),
where  \(\mathsf{I}_\varrho\) is an \emph{induction operator}: this is the inductive counterpart to the notation~\(\iota_xA(x)\) for the~\(x\) such that~\(A(x)\).

On the other hand, although generalised inductive definitions are presented as ``the method of \citealt{lorenzen55}'' by \citet[page~48]{lorenzenmyhill59}, we have not been able to spot a discussion of these in this book. See \citet[§~2]{coquand21} for the generalised inductive definition of the perfect kernel of a closed set provided by \citet{MR0104569} as a constructive analysis of the Cantor--Bendixson theorem.
See also the discussion by \citet[§51, ``Predicativity and inductive definitions'']{parsons08}
.

\paragraph{The law of excluded middle.}
\label{sec:le-tiers-exclu}

Lorenzen considers that the proof of consistency of a calculus secures the use of the law of excluded middle in this calculus:
\begin{quotation}
  In this book we shall use the logical particles according to the
  rules of so-called \emph{classical logic}. [\dots]
  
  Our treatment can be called ``constructive'' even so, since
classical logic can itself be justified from the standpoint of a
\emph{constructive logic} by way of so-called consistency proofs (cf.\
\citealp{lorenzen62}). Furthermore, for indefinite quantifiers we shall
use only constructively valid inferences. \citep[page~15.]{MR0346097}
\end{quotation}
Another way of justifying the recourse to the law of excluded middle is that the proof of consistency provides a constructive understanding of classical truth within the calculus considered, i.e.\ of the definite quantifiers governed by classical logic. 
\begin{quotation}
  [\dots] recourse
to this law is ``innocuous'' in the following sense. By the use of
excluded middle one never gets statements contradictory to those
provable without excluded middle. 
\citep[pages~33--34.]{MR0346097}
\end{quotation}

Lorenzen therefore uses freely the law of excluded middle for natural numbers but not for real numbers.
%
%
%
%
This should be compared with the Limited Principle of Omniscience (LPO) of \citet{B67}.
An analysis of LPO together with dependent choice is made by \citet{coquandpalmgren00} and with more particulars by \citet{feferman12}, and then by \cite{rathjen19}.





\paragraph{Physics.}
\label{sec:physique}

\citet[page~150]{lorenzen86b} emphasises that Weyl wants mathematics to be applicable in physics and that this is decisive in his enthusiasm for operative mathematics.

The journal issue \citet{schlaudt14b} deals with the common grounds of Dingler, Weyl, and Lorenzen in the foundation of physics. We limit ourselves to the following quotation of  \citealp{lorenzen64} in the translation of \citealt[page~235]{lorenzen87b}.



\begin{quotation}
  Of course, we have accomplished nothing by merely asserting the a priori character of protophysics. Like Kant, we must ask what makes an a priori protophysics possible. Unfortunately, the Kantian investigations of this question are completely unsatisfactory in their details. Most physicists and mathematicians
  fail to see any problem here---or even dogmatically reject the possibility of the problem. Hugo Dingler and Hermann Weyl are recent exceptions.
\end{quotation}

\paragraph{Abstraction and its reception by Angelelli.}

Ignacio Angelelli proposes another emphasis of Lorenzen's reception of Weyl's work: abstract objects.

\begin{quotation}
  In the history of philosophy the word ``abstraction'' and cognate expressions (``abstract'',
etc.)\ have had a \emph{genuine} meaning, according to which abstraction involves an
operation by which something is \emph{retained} and something else is \emph{left out} [\dots]. In the special sense relevant for philosophy, the operation is intellectual, and the retaining and leaving out pertain to
our mental consideration of things. \cite[page~11.]{angelelli04}
\end{quotation}

He underlines the great relevance of the conception expounded in \citealt{lorenzen55}.
\begin{quotation}
  Thereby [that is according to what is customary since Frege and Russell] abstrac\bs
tion is to be reduced to the introduction of ``classes''. However, we will see below
that classes are nothing else but a special kind of abstract objects. \citep[Translation of {\citealt[page~101]{lorenzen55}}, in][page~26.]{angelelli04}
\end{quotation}

\citet{lorenzen55} introduces abstraction as follows: to define abstract objects as given concrete objects together with an equivalence relation that makes some of them equal; a statement about the abstract object is a statement about the concrete object that is invariant with respect to the equivalence relation. This should be compared to Bishop's conception of equality in constructive mathematics \citep[see][page~13]{B67}. In particular, classes (i.e.\ sets) of numbers are introduced as the abstraction of arithmetical predicates for the relation of logical equivalence, as does Weyl in our first citation of \emph{The continuum} on page~\pageref{pred}.

\citet[and in all his writings on abstraction cited there]{angelelli04} credits Lorenzen for having set up a modern theory of abstraction. He does so as well in several reviews.
\begin{quotation}
  Thus Weyl's view of sets is neither the naive Cantorian nor the axiomatic; Weyl wants his sets to emerge from properties. If the reader asks ``how?'' no clear answer will be found in this volume, except for hints at ``ideal elements'' and ``definitions by abstraction''. \citet[pp.~45–47]{weyl87} rightly points to Frege as his source in this respect but he is confused about abstraction; Frege did not really do abstraction 
  [\dots]. Here again the consideration of Lorenzen's work is recommended, at any rate for the understanding of what is still unclear or undeveloped in Weyl. (\citealp{angelelli91}.)\medskip
  
  The traditional and Cantorian insistence on number as an abstractum was basically right; it only needed a better theory of abstraction. This was done, apparently for the first time, by P. \citet[Section~13]{lorenzen55}. (\citealp{angelelli86}.)\medskip

The story begins with Peano who, contrary to Frege, was interested in abstraction and gave a place to it in logical theory under the label of ``definitions by abstraction''. 
[\dots] To be sure, Frege's (two-stage!) method can be ``overhauled'', or reconstructed as genuine abstraction. To this end, first, Frege's method should be dismantled and, secondly, rebuilt, according to the guidelines sketched by Peano, praised but only hinted at by Hermann Weyl, and systematically developed by Paul Lorenzen. In the latter's book [\citealp[§6]{MR0209410}], or in his earlier work [\citealp{lorenzen55}], real numbers are in fact introduced by abstraction---by genuine abstraction, that is, not by any pseudo-abstraction. 
(\citealp{angelelli01}.)
\end{quotation}

Lorenzen, as \citet[page~15]{B67}, defines real numbers as Cauchy sequences together with the equivalence relation that their difference is a null sequence, instead of resorting to equivalence classes, i.e.\ quotients. The consideration of equivalence classes leads to the difficulty of choosing representatives and e.g.\ requires countable choice for keeping constructive the proof by \citet[page~25]{B67} that the real numbers are uncountable. Equality in the \citet{hottbook} proposes a way to define quotients that solves this issue.
\[*\ *\ *\]

Hermann Weyl and Paul Lorenzen share a profound interest into physics and especially into the continuum. They consider that mathematical analysis must be built up on the concept of number by a process called by the former ``the mathematical process'' and by the latter ``logical reflection'', leading to the concepts of function, relation, inductive definition. Both call for retrieving the corpus of classical analysis without resorting to impredicative concept formations like Dedekind cuts and do so by deepening the logical analysis of the objects involved, in particular by a better understanding of abstraction. As Poincaré asserted, the real numbers form an indefinite totality; \citet{weyl18} concludes that the least-upper-bound property ``must certainly be abandoned''; \citet{MR0043748,lorenzen55} retrieves it by introducing language levels; \citet{MR0209410} avoids them by considerations of definiteness. 

Weyl's enthusiasm for Lorenzen's operative mathematics witnesses how new logical insights can provide ``pillars of enduring strength'' for analysis and how logic, by specifying the way mathematical objects are given to us, is instrumental in this justification.

\paragraph{Acknowledgment.}

The authors thank Mark van Atten for providing them with the administrative documents relative to Lorenzen's visit to the Institute for Advanced Study.
They thank Jan von Plato for sharing his findings about Lorenzen in the Gödel papers. 

\bibliography{../Lorenzen.bib,lorenzen-weyl.bib,../lorenzen_s_algebraic_and_logistic_investigations_on_free_lattices/lorenzen-algebraic_and_logistic_investigations_on_free_lattices.bib}

\end{document}